# On the theory of minimal submanifolds and harmonic maps from the point of view of the generalized Bochner technique

Sergey E. Stepanov[a,b] *, Irina I. Tsyganok[b]

[a] *Department of Mathematics, Russian Institute for Scientific and Technical Information of the Russian Academy of Sciences, 20, Usievicha street, 125190 Moscow, Russia*

[b] *Department of Mathematics, Finance University, 49-55, Leningradsky Prospect, 125468 Moscow, Russia*

**Abstract:** In the present paper, we study harmonic mappings of complete Riemannian manifolds, as well as minimal and stable minimal submanifolds of complete Riemannian manifolds. We examine classical theorems in the theory of these manifolds from the perspective of the generalized Bochner technique.

> Unfortunately, our knowledge of the subject is not mature enough to give applications to solve outstanding problems in geometry, such as the Hodge conjecture. But the future is bright.
>
> S.-T. Yau

## 1. Introduction

In this paper, we study minimal submanifolds and harmonic maps from the perspective of the generalized Bochner technique, which are long-standing variational problems in Riemannian geometry. The generalized Bochner technique (see, for instance, [1]) is an important component of geometric analysis (see, for example, [2] and [3]). It has many valuable applications in several important areas of Riemannian geometry at large (see, for instance, [4]) and in General Relativity (see, for instance, [5]). Typically, applying the generalized Bochner technique leads to the Yau and Schoen theorem (see [8, p. 337]), which generalizes the Eells and Sampson theorem to the case of non-existence of harmonic mappings with finite Dirichlet energy on complete Riemannian manifolds. In the next section, we will consider this example in more detail and present a new formulation of the Yau and Schoen theorem.

---

* Corresponding author.
E-mail address: s.e.stepanov@mail.ru (S.E. Stepanov)

On the other hand, obtaining rigidity results for minimal and, in particular, stable minimal submanifolds in Riemannian manifolds has been a problem of significant interest in the geometric analysis community. This has been explored from the classical paper [9] to more recent works, such as [10], which also utilized the generalized Bochner technique. As an example, we recall the following theorem (see [11, Theorem 3.1]): Let $(M, g)$ be an $n$-dimensional stable manifold minimally immersed in an $(n + 1)$-dimensional Riemannian manifold $(\bar{M}, \bar{g})$ with non-negative sectional curvature. If $(M, g)$ is parabolic, then it must be totally geodesic and have non-negative sectional curvature. We will continue this research in the fourth section of our article. Bochner's generalized methodology enables us to obtain new results in this actively studied area of research (see, for example, [12]). We believe that Bochner's method will help find applications for solving unsolved problems in geometry and thus pave the way to a brighter future.

We presented some of our results in reports at the "XIV Belarusian Mathematics Conference" (October 28–November 1, 2024, Minsk, Belarus).

## 1. Harmonic mapping of a complete manifold into a manifold with non-negative sectional curvature

A systematic study of harmonic maps was initiated in 1964 by Eells and Sampson [7]. Detailed presentations of the results can be found in [8], [13] - [17], and many other publications. For definitions, notations, and results, we will refer to these works.

Let $(M, g)$ and $(\bar{M}, \bar{g})$ be two connected smooth Riemannian manifolds. The Dirichlet energy, i.e., the stretching energy of a smooth map $f: (M, g) \to (\bar{M}, \bar{g})$, is formally defined as: $E(f) = \frac{1}{2} \int_M e(f)(x) \, dv_g$, where $dv_g$ denotes the volume element of $(M, g)$, and $e(f)(x) := (trace_g (f^{-1}\bar{g}))(x)$ is a non-negative scalar. The term $e(f)(x)$ is known as the *energy density* of $f$ at $x$, and it provides a measure of how much the map $f$ distorts or stretches the metric $\bar{g}$ of the target manifold $\bar{M}$ at each point in $M$. In addition, $f^{-1}\bar{g}$ is a tensor induced on $M$ from $\bar{g}$ by $f$.

Let $D$ be the pull-back covariant derivative on the bundle $T^*M \otimes f^{-1}TN$ over $M$ defined by Levi-Civita connections of $(M,g)$ and $(\bar{M},\bar{g})$, one has a bilinear map $(Ddf)(x)$ of tangent spaces $T_xM \times T_xM \to T_{f(x)}\bar{M}$ at each point $x \in M$. This section is the Hessian of $f$ and has been called *second fundamental form* of $f$ (see [1, p. 2]).

A smooth map $f: (M,g) \to (\bar{M},\bar{g})$ is called harmonic if it is a critical point of the $E(f)$. Furthermore, $f$ is harmonic if and only if $trace_g(Ddf) = 0$ (see the definition and examples of harmonic map [1, pp. 116-118] and [13, pp. 293-294]). Take coordinate neighborhoods $U \subset M$ with local coordinates $x^1, \ldots, x^n$ and $\bar{U} \subset \bar{M}$ with local coordinates $\bar{x}^1, \ldots, \bar{x}^m$ such that $f(U) \subseteq \bar{U}$. We denote by $g_{ij}$ the local components of the Riemannian metric $g$ on $U \subset M$ and by $\bar{g}_{\alpha\beta}$ those of the Riemannian metric $\bar{g}$ on $\bar{U} \subset \bar{M}$ where the indices $i,j,k,l$ run over the range $\{1,\ldots,n\}$ and $\alpha,\beta,\gamma,\varepsilon$ run over the range $\{1,\ldots,m\}$. Suppose that $f: (M,g) \to (\bar{M},\bar{g})$ is given by equations $\bar{x}^\alpha = f^\alpha(x^1, \ldots, x^n)$ with respect to the local coordinates of $U \subset M$ and $\bar{U} \subset \bar{M}$. We put $f_i^\alpha = d\bar{x}^\alpha/dx^i$ then the differential $df$ of the mapping $f$ is represented by the matrix $(f_i^\alpha)$ with respect to the local coordinates of $U \subset M$ and $\bar{U} \subset \bar{M}$. In this case, we can rewrite $f^{-1}\bar{g}$ and $e(f)(x)$ as $(f^{-1}\bar{g})_{kl} = f_k^\alpha f_l^\beta \bar{g}_{\alpha\beta}$ and $e(f)(x) = \frac{1}{2} g^{ij}(f_i^\alpha f_j^\alpha \bar{g}_{\alpha\beta})$, respectively, where $(g^{ij}) = (g_{ij})^{-1}$.

We recall that if $f$ is harmonic, then the *Weitzenböck–Bochner formula*

$$(\Delta e(f))(x) = \|(Ddf)\|_{\tilde{g}}^2(x) + Q(f)(x), \qquad (2.1)$$

holds (see, for example, [7, p. 12]; [1, p. 12]; [14, p. 506] and etc.), where $\Delta e(f)(x) := (trace_g (\text{Hess } e(f)))(x)$ is the *Laplace-Beltrami operator* to the energy density $e(f)(x)$ and $\|(Ddf)\|_{\tilde{g}}^2(x)$ is the square of the norm of the second fundamental form of $f$ with respect to the metric $\tilde{g}$ on $T^*M \otimes T^*M \otimes f^{-1}T\bar{M}$ induced by the metrics $g$ and $\bar{g}$. In turn, $Q(f)(x)$ has the form (see [1, p. 3]; [7, p. 123]; [14, p. 506] and [17, p. 244])

$$Q(f)(x) = -\left(f_i^\alpha f_j^\beta f_k^\gamma f_l^\varepsilon g^{ik}g^{jl}\right)\bar{R}_{\alpha\beta\gamma\varepsilon} + g^{ik}g^{jl}\left(f_k^\alpha f_l^\beta \bar{g}_{\alpha\beta}\right)R_{ij} \qquad (2.2)$$

for the local components $\bar{R}_{\alpha\beta\gamma\varepsilon}$ of the Riemannian curvature tensor $\bar{R}$ of $(\bar{M}, \bar{g})$ and the local components $R_{ij}$ of the Ricci tensor $Ric$ of $(M, g)$ with respect to the local coordinates of $\bar{U} \subset \bar{M}$ and $U \subset M$, respectively. Suppose that the Ricci curvature of $(M, g)$ is non-negative and the sectional curvature of $(\bar{M}, \bar{g})$ is non-positive, then from (2.2) we obtain $Q(f)(x) \geq 0$ (see, for instance, [7, p. 123]; [8, p. 336]; [17, p. 246]). Furthermore, we recall the celebrated vanishing theorem on harmonic maps from [7, p. 124] which state the following: if $f:(M, g) \to (\bar{M}, \bar{g})$ is any harmonic mapping between a compact Riemannian manifold $(M, g)$ with non-negative Ricci tensor and a Riemannian manifold $(\bar{M}, \bar{g})$ with non-positive sectional curvature, then $f$ is *totally geodesic* (see [7, p. 9]) and has constant the pointwise energy density $e(f)(x)$ (see [7, p. 10]; [1, p. 1]). Furthermore, if there is at least one point of $(M, g)$ at which its Ricci curvature is positive, then every harmonic map $f:(M, g) \to (\bar{M}, \bar{g})$ is constant (see also [7]). The most recent vanishing theorem proved in article [8, p. 337]. Let us formulate this theorem.

**Theorem 2.1**. *Let $(M, g)$ be a complete Riemannian manifold with Ricci curvature $Ric \geq 0$ at each point of $M$ and $(\bar{M}, \bar{g})$ be a compact Riemannian manifold with sectional curvature $\overline{sec} \leq 0$ at each point of $\bar{M}$, then any harmonic map $f:(M, g) \to (\bar{M}, \bar{g})$ with finite Dirichlet energy $E(f)$ must be homotopic to a constant map.*

Now, for a map $f:(M, g) \to (\bar{M}, \bar{g})$ between a Riemannian manifolds, we denote by $Ric(X)$ the Ricci curvature of $(M, g)$ for the unit tangent vector $X \in T_x M$ at $x \in M$. In turn, we denote by $\overline{Ric}(\bar{X})$ the Ricci curvature of $(\bar{M}, \bar{g})$ for the unit tangent vector $\bar{X} \in T_{\bar{x}}\bar{M}$ at $\bar{x} \in \bar{M}$ and by $\overline{sec}(\bar{\pi})$ the sectional curvature of $(M, g)$ for the tangent 2-plane $\bar{\pi} \subset T_{\bar{x}}\bar{M}$ at $\bar{x} \in \bar{M}$. Next, set

$$Ric_{min}(x) := min_{X \in T_x M} Ric(X), \quad \overline{Ric}_{max}(x) := max_{\bar{X} \in T_{\bar{x}}\bar{M}} \overline{Ric}(\bar{X}),$$

$$\overline{sec}_{min}(x) := min_{\bar{\pi} \subset T_{\bar{x}}\bar{M}} \overline{sec}(\bar{\pi}).$$

Then, in contrast to the previous theorem, we prove our vanishing theorem for harmonic maps in a new formulation.

**Theorem 2.2.** *Let $f:(M, g) \to (\bar{M}, \bar{g})$ be a harmonic map between complete and compact Riemannian manifolds, respectively, with finite Dirichlet energy $E(f) =$*

$\frac{1}{2}\int_M e(f)(x)\, dv_g$, where $e(f)(x)$ denotes the pointwise energy density of $f$. Assume the following condition: For all $x \in M$ the inequalities

$$\overline{sec}_{min}(x) \cdot e(f)(x) \leq Ric_{min}(x), \qquad \overline{sec}_{min}(x) \geq \frac{1}{m}\overline{Ric}_{max}(x) \quad (2.3)$$

hold, where $Ric_{min}(x)$ denotes the minimal Ricci curvature of $(M,g)$ at $x$, $e(f)(x)$ denotes the pointwise energy density of $f$ at an arbitrary point $x \in M$, $\overline{Ric}_{max}(x)$ and $\overline{sec}_{min}(x)$ denote the maximal Ricci and minimal sectional curvatures of $(\bar{M},\bar{g})$ at $f(x)$, respectively. Then the harmonic map $f$ is homotopic to constant on each compact set in $M$.

*Proof.* We set $\Phi^{\alpha\beta} := f_k^\alpha f_l^\beta g^{kl}$ and rewrite (2.2) in the form

$$Q(f)(x) = -\bar{R}_{\alpha\beta\gamma\varepsilon}\Phi^{\alpha\gamma}\Phi^{\beta\varepsilon} + g^{ik}g^{jl}\left(f_k^\alpha f_l^\beta \bar{g}_{\alpha\beta}\right)R_{ij}$$

with respect to the local coordinates of $\bar{U} \subset \bar{M}$ and $U \subset M$, respectively, where. Then we diagonalize the symmetric tensor $\Phi$ with local components $\Phi^{\alpha\beta}$, using an orthonormal basis $\{\bar{e}_1, \ldots, \bar{e}_m\}$ at each point $f(x) \in \bar{M}$. Then by [18, p. 436] and [19, p. 391] we have

$$\bar{R}_{\alpha\beta}\bar{g}_{\gamma\varepsilon}\Phi^{\alpha\gamma}\Phi^{\beta\varepsilon} - \bar{R}_{\alpha\beta\gamma\varepsilon}\Phi^{\alpha\gamma}\Phi^{\beta\varepsilon} = \sum_{\alpha<\beta}\overline{sec}\,(\bar{e}_\alpha,\bar{e}_\beta)(\bar{\lambda}_\alpha - \bar{\lambda}_\beta)^2 \geq$$

$$\geq \overline{sec}_{min}(x)\sum_{\alpha<\beta}^{1,2,\ldots,m}(\bar{\lambda}_\alpha - \bar{\lambda}_\beta)^2 = \overline{sec}_{min}(x)\left(m\sum_{\alpha=1}^m(\bar{\lambda}_\alpha)^2 - \left(\sum_{\alpha=1}^m \bar{\lambda}_\alpha\right)^2\right),$$

where $\bar{\lambda}_1, \ldots, \bar{\lambda}_m$ are eigenvalue of $\Phi$ and $\overline{sec}\,(\bar{e}_\alpha,\bar{e}_\beta)$ is the sectional curvature in the direction of $\bar{\pi}(x) = span\{\bar{e}_\alpha,\bar{e}_\beta\} \subset T_{f(x)}\bar{M}$. In this case, a straightforward calculation using (2.2) yields

$$Q(f)(x) = \sum_{\alpha<\beta}\overline{sec}\,(\bar{e}_\alpha,\bar{e}_\beta)(\bar{\lambda}_\alpha - \bar{\lambda}_\beta)^2 + \bar{g}_{\alpha\beta}f_k^\alpha f_l^\beta g^{ik}g^{jl}\left(R_{ij} - \bar{R}_{\gamma\varepsilon}f_i^\gamma f_j^\varepsilon\right) \geq$$

$$\geq \overline{sec}_{min}(x)\left(m\sum_{\alpha=1}^m(\bar{\lambda}_\alpha)^2 - \left(\sum_{\alpha=1}^m \bar{\lambda}_\alpha\right)^2\right) +$$

$$+ \left(\bar{g}_{\alpha\beta}f_k^\alpha f_l^\beta\right)g^{ik}g^{jl}\left(Ric_{min}(x)\cdot g_{ij} - \overline{Ric}_{max}(x)\cdot\left(\bar{g}_{\gamma\varepsilon}f_i^\gamma f_j^\varepsilon\right)\right) =$$

$$\geq \overline{sec}_{min}(x)\left(m\sum_{\alpha=1}^m(\bar{\lambda}_\alpha)^2 - \left(\sum_{\alpha=1}^m \bar{\lambda}_\alpha\right)^2\right) +$$

$$+ Ric_{min}(x)\bar{g}_{\alpha\beta}\Phi^{\alpha\beta} - \overline{Ric}_{max}(x)\bar{g}_{\alpha\beta}\Phi^{\alpha\gamma}\bar{g}_{\gamma\varepsilon}\Phi^{\beta\varepsilon} =$$

$$= \overline{sec}_{min}(x)\left(m\sum_{\alpha=1}^m(\bar{\lambda}_\alpha)^2 - \left(\sum_{\alpha=1}^m \bar{\lambda}_\alpha\right)^2\right) +$$

$$+ Ric_{min}(x) \left(\sum_{\alpha=1}^{m} \bar{\lambda}_\alpha\right) - \overline{Ric}_{max}(x) \sum_{\alpha=1}^{m} (\bar{\lambda}_\alpha)^2 =$$

$$= \left(\sum_{\alpha=1}^{m} (\bar{\lambda}_\alpha)^2\right) \left(m \, \overline{sec}_{min}(x) - \overline{Ric}_{max}(x)\right) +$$

$$+ \left(\sum_{\alpha=1}^{m} \bar{\lambda}_\alpha\right) \left(Ric_{min}(x) - \overline{sec}_{min}(x) \left(\sum_{\alpha=1}^{m} \bar{\lambda}_\alpha\right)\right),$$

where

$$\sum_{\alpha=1}^{m} \bar{\lambda}_\alpha = trace_{\bar{g}} \Phi = \left(f_k^\alpha f_l^\beta g^{kl}\right) \bar{g}_{\alpha\beta} = \left(f_k^\alpha f_l^\beta \bar{g}_{\alpha\beta}\right) g^{kl} = e(f)(x),$$

$$\sum_{\alpha=1}^{m} (\bar{\lambda}_\alpha)^2 = trace_{\bar{g}} \Phi^2 = \Phi^{\alpha\beta} \Phi^{\gamma\delta} \bar{g}_{\alpha\gamma} \bar{g}_{\beta\delta} =$$

$$= \left(f_i^\alpha f_j^\beta g^{ij}\right) \left(f_k^\gamma f_l^\delta g^{kl}\right) \bar{g}_{\alpha\gamma} \bar{g}_{\beta\delta} =$$

$$= \left(f_i^\alpha f_k^\gamma \bar{g}_{\alpha\gamma}\right) \left(f_j^\beta f_l^\delta \bar{g}_{\beta\delta}\right) g^{ij} g^{kl} = \|(f^{-1}\bar{g})(x)\|_g^2(x).$$

Therefore, if the pointwise inequalities

$$\overline{sec}_{min}(x) \cdot e(f)(x) \leq Ric_{min}(x), \qquad \overline{sec}_{min}(x) \geq \frac{1}{m} \overline{Ric}_{max}(x)$$

hold at each point $x \in M$, then $Q(f)(x) \geq 0$. Therefore, from (2.1) we obtain $\Delta e(f)(x) \geq 0$ at each point $x \in M$ and hence $e(f)(x)$ is a *subharmonic function*. In the paper [20], Yau has shown that every non-negative $L^2$-integrable subharmonic function on a complete Riemannian manifold must be constant. Applying this to $\sqrt{e(f)(x)}$, we conclude that $e(f)(x)$ is a constant function (see also [8, p. 337]). On the other hand, every complete non-compact Riemannian manifold with nonnegative Ricci curvature has infinite volume (see [30]). In our case, we have $Ric \geq 0$ since $Ric_{min}(x) \geq e(f)(x) \cdot \overline{sec}_{min}(x) \geq 0$ at each point $x \in M$ (see our Remark 2.1). Therefore, the volume of $(M, g)$ is infinite. This forces the constant function $e(f)(x)$ to be zero and hence $f$ to be a constant map (see also [8, p. 337]). In conclusion, we refer to Theorem 1 of Yau and Shoen and can conclude that our theorem is also true.

**Remark 2.1.** In the conditions of our theorem, we considered the following inequality: $\overline{Ric}_{max} \leq m \, \overline{sec}_{min}$ at an arbitrary point $x \in M$. Since the inequality $\overline{Ric}(\bar{X}) \geq (m-1) \overline{sec}_{min}(x)$ holds for an arbitrary unit vector $X \in T_x M$ at each point $x \in M$, we obtain the double inequality:

$$(m-1) \overline{sec}_{min}(x) \leq \overline{Ric}(\bar{X}) \leq m \, \overline{sec}_{min}(x).$$

This double inequality implies that $\overline{sec}_{min}(x) \geq 0$ at any point $f(x) \in \overline{M}$. Therefore, $(\overline{M}, \bar{g})$ is a Riemannian manifold of positive sectional curvature at each point of $f(M) \subset \overline{M}$ in our theorem (compare to the theorem of Yau and Shoen).

**Remark 2.2.** From inequalities (2.3) we deduce the condition

$$Ric_{min}(x) \geq \frac{1}{m} e(f)(x) \overline{Ric}_{max} \geq 0$$

at each point $x \in M$. This inequality can be considered as a generalization of other inequality $Ric \geq f^{-1}\overline{Ric}$ obtained in our article [16] as a condition for the obstruction of the existence of a harmonic map $f: (M, g) \to (\overline{M}, \bar{g})$.

Additionally, let $(M, g)$ be a compact Riemannian manifold, then we can apply Stokes' theorem $\int_M \Delta e(f)(x) \, dv_g = 0$ to (2.1). As a result, we obtain the integral formula

$$\int_M \left( \|(Ddf)\|_{\bar{g}}^2 (x) + Q(f)(x) \right) dv_g = 0.$$

and hence the following inequality

$$\int_M Q(f)(x) \, dv_g \leq 0. \tag{2.4}$$

In this case, if the inequalities

$$\overline{sec}_{min}(x) \cdot e(f)(x) \leq Ric_{min}(x), \qquad \overline{sec}_{min}(x) > \frac{1}{m} \overline{Ric}_{max}(x)$$

hold at each point $x \in M$, then $Q(f)(x) \geq 0$. Then we can conclude from inequality (2.4) that $Q(f)(x) = 0$ at each point $x \in M$. Therefore, in particular, $\|(f^{-1}\bar{g})(x)\|_{\bar{g}}^2(x) = 0$ at each point $x \in M$. Hence $f$ is a constant map. Moreover, we recall (see [19]) that an $m$-dimensional ($m \geq 3$) connected compact Riemannian manifold $(\overline{M}, \bar{g})$, whose Ricci tensor and sectional curvature satisfy the strict inequality $\overline{sec}_{min}(x) > \frac{1}{m}\overline{Ric}_{max}(x)$ at each point of $(\overline{M}, \bar{g})$, is diffeomorphic to the spherical space form $\mathbb{S}^m / \Gamma$. Furthermore, if $(\overline{M}, \bar{g})$ is simply connected, then it is diffeomorphic to $\mathbb{S}^m$. As a result, the new version of the Eells and Sampson theorem on the vanishing of harmonic maps holds.

**Corollary 2.3.** *Let $f: (M, g) \to (\overline{M}, \bar{g})$ be a harmonic map between compact Riemannian manifolds $(M, g)$ and $(\overline{M}, \bar{g})$. If the inequalities*

$$Ric_{min}(x) \geq \overline{sec}_{min}(x) \cdot e(f)(x), \qquad \overline{sec}_{min}(x) > \frac{1}{m} \overline{Ric}_{max}(x) \tag{2.3}$$

hold at each point $x \in M$, where $Ric_{min}(x)$ denotes the minimal Ricci curvature of $(M, g)$ at $x$, $e(f)(x)$ denotes the pointwise energy density of $f$ at an arbitrary point $x \in M$, $\overline{Ric}_{max}(x)$ and $\overline{sec}_{min}(x)$ denote the maximal Ricci and minimal sectional curvatures of $(\bar{M}, \bar{g})$ at $f(x)$, respectively. Then $(M, g)$ is diffeomorphic to the spherical space form $\mathbb{S}^m/\Gamma$, and the harmonic map $f$ is constant. Furthermore, if $(\bar{M}, \bar{g})$ is simply connected, then it is diffeomorphic to $\mathbb{S}^m$.

## 3. Complete minimal submanifolds in Riemannian manifolds

We recall that one of the most intriguing topics in the calculus of variations within Riemannian geometry is the study of minimal submanifolds $(M, g)$ of a Riemannian manifold $(\bar{M}, \bar{g})$ (see, for instance, [21] and [22, § 5.1]). To be precise and to establish our notation, we note the following. Let $(M, g)$ be an $n$-dimensional complete manifold isometrically immersed in an $(n + k)$-dimensional Riemannian manifold $(\bar{M}, \bar{g})$ of constant curvature $C$ for $k \geq 1$. We denote by $\nabla$ and $\bar{\nabla}$ the Levi-Civita connections on $(M, g)$ and $(\bar{M}, \bar{g})$, $(\bar{M}, \bar{g})$, respectively.

For vector fields $X, Y \in C^\infty(TM)$, the *Gauss formula* of $(M, g) \subset (\bar{M}, \bar{g})$ has the form $\varphi(X, Y) = \bar{\nabla}_X Y - \nabla_X Y$, where $\varphi$ is a normal-bundle symmetric 2-tensor field on $M$, which is called the second fundamental form of $(M, g)$. In this case $H(x) = \frac{1}{n}\sum_{i=1}^{n} \varphi(e_i, e_i)$ is the mean curvature vector, where $e_1, \ldots, e_n$ is an orthonormal basis of $T_x M$ at each point $x \in M$. (see, for instance, [22, §3.3]; [21, p. 68]). Then $(M, g)$ is called a minimal submanifold of $(\bar{M}, \bar{g})$ if the mean curvature vector field vanishes identically (see also [21, p. 68] and [22, §3.3]). We denote by $\tilde{g}$ a metric of the bundle $TM \oplus T^\perp M$ over $M$, then, in particular, the square of the length of the second fundamental form $\varphi$ of $(M, g)$ has the form

$$\|\varphi\|_{\tilde{g}}^2(x) = \sum_{i=1}^{n}\sum_{j=1}^{n} \bar{g}\left(\varphi(e_i, e_j), \varphi(e_i, e_j)\right),$$

where $e_1, \ldots, e_n$ is an orthonormal basis of $T_x M$ at each point $x \in M$. It easy that $\|\varphi\|_{\tilde{g}}^2(x)$ is a smooth function defined on $M$. In turn, a hypersurface $(M, g)$ is called *totally umbilical* if $\varphi = H(x)g$ at each point $x \in M$. On the other hand, a hypersurface $(M, g)$ is a *totally geodesic* submanifold of $(\bar{M}, \bar{g})$ if its second

fundamental form is vanishing. Let us recall the following classical statement due to Chern, do Carmo and Kobayashi (see [9; Theorem 1 and its Corollary]).

**Theorem 3.1.** *Let $(M, g)$ be an $n$-dimensional compact Riemannian manifold that is minimally immersed in an $(n + k)$-dimensional Riemannian manifold $(\bar{M}, \bar{g})$ with constant curvature $C$. If $(M, g)$ is not totally geodesic and $\|\varphi\|_{\bar{g}}^2(x) \leq \frac{kn}{2k-1} C$, where $\varphi$ represents the second fundamental form of $(M, g)$, then $(M, g)$ is parallel submanifold and $\|\varphi\|_{\bar{g}}^2(x) = \frac{kn}{2k-1} C$.*

**Remark 3.1.** Recall that a submanifold $(M, g)$ of $(\bar{M}, \bar{g})$ is said to be *parallel* if its second fundamental form $\varphi$ is parallel with respect to the *Van der Waerden - Bortolotti connection* $\widetilde{\nabla}$ defined on the bundle $TM \oplus T^\perp M$ (see details in [22, Chapter 8] and [23]).

Let us consider further complete Riemannian manifolds. It is well known that on a complete noncompact Riemannian manifold, there always exists a Green's function (see [24, Paragraph 4.2]). The existence and nonexistence of a positive Green's function divide the class of complete manifolds into two categories. In general, the methods used in dealing with function theory on these manifolds are different. Namely, a Riemannian manifold is called parabolic if it admits no positive Green's function, and hyperbolic otherwise (see [2, Definition 20.1]). An example of a parabolic manifold is a complete Riemannian manifold $(M, g)$ of finite volume (see [25]).

**Remark 3.2.** In contrast to the above, Grigoryan does not use the adjective "hyperbolic" as an antonym to "parabolic", since, as he said, there are other generalizations of the concept of hyperbolicity, for example, Gromov's hyperbolicity (see [24, p. 165]).

Using the information on parabolic manifolds, we will show the following generalized theorem of Chern, Kobayashi, and Do Carmo.

**Theorem 3.2.** *Let $(M, g)$ be an $n$-dimensional complete parabolic manifold (in particular, a complete Riemannian manifold of finite volume) that is minimally immersed in an $(n + k)$-dimensional Riemannian manifold $(\bar{M}, \bar{g})$ with constant*

*curvature $C$. If $(M,g)$ is not totally geodesic and $\|\varphi\|_{\bar{g}}^2(x) \leq \frac{kn}{2k-1} C$, where $\varphi$ represents the second fundamental form of $(M,g)$, then $(M,g)$ is parallel submanifold and $\|\varphi\|_{\bar{g}}^2(x) = \frac{kn}{2k-1} C$.*

*Proof.* Let $(M,g)$ be an $n$-dimensional Riemannian manifold which is isometrically immersed in an $(n+q)$-dimensional Riemannian manifold $(\bar{M}, \bar{g})$ with constant curvature $C$. According to Eq. (3.10) and Eq. (3.12) from [9, pp. 65-66], we have the Weitzenböck–Bochner formula

$$\tfrac{1}{2} (\Delta \|\varphi\|_{\bar{g}}^2)(x) = \|\tilde{\nabla}\varphi\|_{\bar{g}}^2(x) + \left(n\, C - \left(2 - \tfrac{1}{k}\right) \|\varphi\|_{\bar{g}}^2(x)\right) \|\varphi\|_{\bar{g}}^2(x), \quad (3.1)$$

where $\tilde{\nabla}\varphi$ is a covariant derivative with respect to the connection $\tilde{\nabla}$ defined on $TM \oplus T^\perp M$ (see [7, §3.1]; [8, p. 62; 66]). Therefore, if we suppose that $(M,g)$ is not totally geodesic and $\|\varphi\|_{\bar{g}}^2(x) \leq n\, C / \left(2 - \tfrac{1}{k}\right)$, then from (3.1) we conclude that $(\Delta \|\varphi\|_{\bar{g}}^2)(x) \geq 0$. In this case, $\|\varphi\|_{\bar{g}}^2(x)$ is a non-negative subharmonic function bounded above. Note that the definition of parabolicity allows the characterization: a complete manifold $(M,g)$ is parabolic if every bounded subharmonic function on $(M,g)$ is a constant (see [24, p. 164]). Therefore, in particular, if $(M,g)$ has finite volume, then $\|\varphi\|(x) = const$. In this case, from (3.1) we obtain $\|\varphi\|_{\bar{g}}^2(x) = n\, C$ and $\tilde{\nabla}\varphi = 0$. In conclusion, we recall that a Riemannian submanifold is said to be parallel if its second fundamental form is parallel with respect to $\tilde{\nabla}$ (see details in [23, Chapter 8]). The proof is complete.

We consider next a smooth connected $n$-dimensional Riemannian manifold $(M,g)$ that is isometrically immersed into an $(n+1)$-dimensional Riemannian manifold $(\bar{M}, \bar{g})$. Such a manifold is called a hypersurface. In this case, the Gauss formula says that $\nabla_X Y = \bar{\nabla}_X Y - g(A_g X, Y)\, N$, where $X, Y$ are any vector fields tangent to $(M,g)$, the vector field $N$ is the global unite vector field normal to $(M,g)$ and $A_g$ stands the shape operator of $(M,g)$. The second fundamental form of a hypersurface $(M,g)$ is defined by the identity $\varphi(X,Y) := g(A_g X, Y)$, where $X, Y$ are any vector fields tangent to $(M,g)$. Then $\varphi$ is a symmetric differential 2-form on $(M,g)$. In this case we denote by $\|\varphi\|_g^2(x) = \sum_{i=1}^n \sum_{j=1}^n \varphi(e_i, e_j) \cdot \varphi(e_i, e_j)$ the square of the length of the second fundamental form of $(M, g)$, where $e_1, \ldots, e_n$ is an orthonormal of $T_x M$ at an arbitrary point $x \in M$. Then the subsequent consequence is obvious.

**Corollary 3.3.** *Let $(M, g)$ be an $n$-dimensional complete parabolic manifold (in particular, a complete manifold with finite volume) that is minimally immersed in an $(n + 1)$-dimensional Riemannian manifold $(\bar{M}, \bar{g})$ with constant curvature $C$. If $(M, g)$ is not totally geodesic and the inequality $\|\varphi\|_g^2(x) \leq n\,C$ holds for its second fundamental form $\varphi$, then $\|\varphi\|_g^2(x) = n\,C$ and $\varphi$ parallel with respect to the Levi-Civita connection $\nabla$ of $(M, g)$.*

By the well-known Eisenhart's theorem, if $\varphi$ is a parallel symmetric differential 2-form on $(M, g)$, then for each point $x \in M$, there is a neighborhood $U \subseteq M$ such that $\varphi = \lambda_1 g_1 + \cdots + \lambda_r g_r$, and locally $(M, g)$ admits a Riemannian direct product structure $(M, g) \supseteq (U, g|U) = (U_1, g_1) \times \ldots \times (U_r, g_r)$, where $(U_1, g_1), \ldots, (U_r, g_r)$ are Riemannian manifolds of dimensions $n_1 \geq 1, \ldots, n_r \geq 1$, respectively, with $n_1 + \cdots + n_r = n$ for some $r$ and $\lambda_1, \ldots, \lambda_r$ are constants (see [26, p. 303]). In this case the conditions $trace_g \varphi = 0$ and $\|\varphi\|_g^2 = n\,C$ can be rewritten as the following

$$n_1 \lambda_1 + \cdots + n_r \lambda_r = 0$$

and

$$n_1 (\lambda_1)^2 + \cdots + n_r (\lambda_r)^2 = n\,C,$$

respectively. Furthermore, according to [9, pp. 67-68], if $(M, g)$ is a minimal hypersurface in a Riemannian manifold $(\bar{M}, \bar{g})$ of constant curvature $C = 1$ and hence $\|\varphi\|_g^2 = n$, then $r = 2$. In this case, we deduce from the above identities that $\lambda_1 = \sqrt{n_2/n_1}$ and $\lambda_2 = -\sqrt{n/n_2}$ for $\dim U_1 = n_1 \geq 1$ and $\dim U_2 = n_2 \geq 1$. Therefore, we have

**Corollary 3.4.** *Let $(M, g)$ be an $n$-dimensional complete parabolic manifold (in particular, a complete manifold with finite volume) that is minimally immersed in an $(n + 1)$-dimensional Riemannian manifold $(\bar{M}, \bar{g})$ with constant curvature 1, then for each point $x \in M$ there is some neighborhood $U \subseteq M$ such that $(M, g)$ is locally a Riemannian direct product of Riemannian manifolds $(U_1, g_1)$ and*

$(U_2, g_2)$ *of constant curvatures and dimensions* $n_1 \geq 1$ *and* $n_2 \geq 1$, *respectively, with* $n_1 + n_2 = n$. *Moreover, in this case* $\varphi = \lambda_1 g_1 + \lambda_2 g_2$ *with* $\lambda_1 = \sqrt{n_2/n_1}$ *and* $\lambda_2 = -\sqrt{n/n_2}$.

If $(\bar{M}, \bar{g})$ is a Riemannian manifold of constant sectional curvature, then the second fundamental form $\varphi$ of $(M, g)$ satisfies the *Codazzi equations* (see [18, p. 436]

$$(\nabla_X \varphi)(Y, Z) = (\nabla_Y \varphi)(X, Z), \qquad (3.2)$$

where $X, Y$ and $Z$ are any vector fields tangent to $(M, g)$. In particular, from Codazzi equations (3.2) we deduce the following equation

$$\delta \varphi = -n \cdot d H, \qquad (3.3)$$

where $\delta \varphi = -\, div\, \varphi$. At the same time, any $\varphi \in C^\infty(S^2 M)$ satisfying Codazzi equations (3.2) is called the *Codazzi tensor* (see [18, p. 435]; [27, p. 350]).

Moreover, if the mean curvature of $(M, g)$ is constant (but not necessarily zero), then the second fundamental form of $(M, g))$ becomes a harmonic bilinear form (see [27, p. 350]), since it satisfies the Codazzi equations and has zero divergence according to (3.3). Using this fact, as well as our vanishing theorem on harmonic bilinear forms on complete manifolds from [28, Theorem 3], let $(M, g)$ be a connected complete noncompact Riemannian manifold with nonnegative sectional curvature. Then there is no nonzero harmonic symmetric bilinear form $\varphi$ that satisfies the condition $\|\varphi\| \in L^p(M)$, which means $\int_M \|\varphi\|_g^p (x)\, dv_g < \infty$ at least for one $p \geq 1$. We arrive at the following conclusion.

**Corollary 3.5.** *Let* $(M, g)$ *be a complete non-compact hypersurface in the Riemannian manifold* $(\bar{M}, \bar{g})$ *with constant sectional curvature. If the sectional curvature of* $(M, g)$ *is non-negative, its mean curvature is constant* (*but not necessary zero*) *and its second fundamental form* $\varphi$ *satisfies the condition* $\|\varphi\|(x) \in L^p(M)$ *for at least one* $p \geq 1$, *then* $(M, g)$ *must be totally geodesic submanifold,*

**Remark 3.3.** In contrast to the above theorem, we recall the following theorem (see [11, Theorem 1.1]): let $(M, g)$ be a complete minimal hypersurface in $(n + 1)$−dimensional Riemannian manifold $(\bar{M}, \bar{g})$ with nonnegative Ricci curvature. If $(M, g)$ is parabolic, then $(M, g)$ must be totally geodesic submanifold of $(\bar{M}, \bar{g})$.

In conclusion of the section, we return to the study of compact minimal submanifolds. It is well-known that for any $n$-dimensional ($n \geq 3$) compact (without boundary) Riemannian manifold $(M, g)$, the algebraic sum $(\operatorname{Im} \delta^* + C^\infty M \cdot g)$ is closed in $S^2 M$, where $\delta^* \theta := \frac{1}{2} L_\xi g$ for the Lie derivative $L_\xi g$ with respect to a smooth vector field $\xi$ which is dual (by $g$) to the 1-form $\theta$ (see [15, p. 35]). In this case, we have the decomposition (see [18, p. 130])

$$S^2 M = (\operatorname{Im} \delta^* + C^\infty M \cdot g) \oplus \left(\delta^{-1}(0) \cap \operatorname{trace}_g^{-1}(0)\right) \quad (3.4)$$

where both factors are infinite dimensional and orthogonal to each other with respect to the $L^2$ inner scalar product $\langle \cdot, \cdot \rangle = \int_M g(\cdot, \cdot) \, dv_g$ for the canonical measure $dv_g$ of $(M, g)$. It's obvious that the second factor $\delta^{-1}(0) \cap \operatorname{trace}_g^{-1}(0)$ of (3.4) is the space of *TT*-tensors.

**Remark 3.4.** Here we recall that a symmetric, divergence-free, and traceless covariant two-tensor is called a *TT*-tensor. As a consequence of a result by Bourguignon, Ebin, and Marsden (see [18, p. 132]), the space of *TT*-tensors is an infinite-dimensional vector space for any closed Riemannian manifold $(M, g)$. Such tensors are of fundamental importance in stability analysis in General Relativity (see, for instance, [29] and [30]) and in Riemannian geometry (see, for instance, [18, p. 346-347]).

Considering the above, we conclude that the second fundamental form of $\varphi$ has the following $L^2$-orthogonal decomposition (see formula (3.3))

$$\varphi = \left(\frac{1}{2} L_\xi g + \lambda g\right) + \varphi^{TT} \quad (3.5)$$

for some vector field $\xi \in C^\infty(TM)$, some *TT*-tensor $\varphi^{TT} \in C^\infty(S^2 M)$ and some scalar function $\lambda(x) \in C^\infty(M)$. In this case, if $\mathring{\varphi}$ denotes the traceless part of $\varphi$, i.e., $\mathring{\varphi} = \varphi - H(x) g$. Then using (3.5) we deduce the following equality

$$\mathring{\varphi} = \varphi + \frac{1}{n}(\delta \theta - n\lambda)g = \left(\frac{1}{2} L_\xi g + \frac{1}{n} \delta \theta \, g\right) + \varphi^{TT}$$

that can be rewritten in the form

$$\mathring{\varphi} = 2S\theta + \varphi^{TT}, \quad (3.6)$$

where $S\theta = L_\xi g + 2/n \, \delta \theta \, g$ is the *Cauchy-Ahlfors operator*. Next, applying $\delta$ to both sides of (3.6), we obtain

$$\delta \overset{\circ}{\varphi} = S^*S \, \theta, \qquad (3.7)$$

for the *Ahlfors Laplacian* $S^*S$ for $S^* := 2\delta$ (see details in [31]). Using (3.3), equation (3.7) can be rewritten in the form

$$\delta \overset{\circ}{\varphi} = -(n-1) \, dH. \qquad (3.8)$$

From (3.1) and (3.8) we deduce the following integral formula

$$\langle S\theta, S\theta \rangle = -(n-1) \int_M (L_\xi H)(x) \, dvol_g. \qquad (3.9)$$

By analyzing formulas (3.6) and (3.9) we can prove our next result, which, in particular, generalizes classical theorem from [21, Theorem 5.4.2].

**Theorem 3.6.** *Let $(M, g)$ be a compact hypersurface in a Riemannian manifold $(\bar{M}, \bar{g})$ with constant sectional curvature and let its second fundamental form $\varphi$ has the $L^2$-orthogonal decomposition*

$$\varphi = \left(\frac{1}{2} L_\xi g + \lambda \, g\right) + \varphi^{TT}$$

*for some vector field $\xi \in C^\infty(TM)$, some TT-tensor $\varphi^{TT} \in C^\infty(S^2M)$ and some scalar function $\lambda(x) \in C^\infty(M)$. If the mean curvature $H(x)$ of $(M, g)$ satisfies the integral equality $\int_M (L_\xi H)(x) \, dvol_g = 0$, then the $L^2$-orthogonal decomposition of $\varphi$ has the form $\varphi = H(x) \, g + \varphi^{TT}$, where $H(x)$ is constant (but not necessary zero), and, in particular, $\overset{\circ}{\varphi}$ is a TT-tensor Codazzi.*

It is well-known the following theorem (see [18, p. 436]): Every Codazzi tensor $\varphi$ with constant trace on a compact Riemannian manifold $(M, g)$ with non-negative sectional curvature is parallel. Furthermore, if the sectional curvatures of $(M, g)$ are positive at some point, then $\varphi$ is a constant multiple of $g$. Therefore, if $(M, g)$ is a hypersurface with *quasi-positive sectional curvature* in a Riemannian manifold $(\bar{M}, \bar{g})$ with constant sectional curvature, then it is a totally umbilical submanifold with constant mean curvature.

A connected Riemannian manifold with constant sectional curvature is called *space form* (see details in [32]). Therefore, in our case, $(M, g)$ must be a spherical space form since it has positive constant curvature. In chapters 4 – 7 of famous monograph [32], Wolf classified spherical space forms as connected Riemannian manifolds locally isometric to the $n$-sphere $\mathbb{S}^n$. In particular, if $(M, g)$ a simply connected manifold, then $M = \mathbb{S}^n$ (see also [27, pp. 200-201]).

Using this fact, along with the results from Theorem 3.6, we arrive at the corollary.

**Corollary 3.7.** *Let $(\bar{M}, \bar{g})$ be an $(n + 1)$-dimensional, where $n \geq 3$, Riemannian manifold with constant sectional curvature and let $(M, g)$ be a compact hypersurface $M \subset \bar{M}$ with constant mean curvature. If $(M, g)$ has quasi-positive sectional curvature, then $(M, g)$ is a spherical space form. Furthermore, if $(M, g)$ is a simply connected manifold, then it is the Euclidean sphere $\mathbb{S}^n$.*

## 4. Complete stable minimal submanifolds in Riemannian manifolds

As we recalled before, stability is related to the non-negativity of the second variation or, equivalently, the non-positivity of the *stability operator* (also called the *Jacobi operator*) on $(M, g)$ (see [10])

$$\mathcal{L} := \Delta + \|\varphi\|_g^2(x) + \overline{Ric}(N, N),$$

where $N$ is a globally defined unit normal vector field on $M$, which means $N_x \in T_x^\perp M$ and $\bar{g}(N_x, N_x) = 1$ at each point $x \in M$. Therefore, in this case we have

$$(\Delta u)(x) \leq - \left( \|\varphi\|_g^2(x) + \overline{Ric}(N_x, N_x) \right) u(x) \qquad (4.1)$$

for any $u(x) \in C^\infty(M)$. Then, using (4.1), we have

$$\tfrac{1}{2}(\Delta u^2)(x) = \|du\|_g^2(x) + u(x)\,(\Delta u)(x) =$$

$$= \|du\|_g^2(x) - \left( \|\varphi\|_g^2(x) + \overline{Ric}(N_x, N_x) \right) u(x)^2 \leq$$

$$\leq - \left( \|\varphi\|_g^2(x) + \overline{Ric}(N_x, N_x) \right) u(x)^2 \qquad (4.2)$$

Therefore, if $u(x)$ is a function that has not zeros on $M$ and satisfying the inequality $\mathcal{L}\,u(x) \leq 0$, and $\overline{Ric}(N_x, N_x) \geq 0$ at every point $x \in M$, then $u^2$ is positive superharmonic function (see [24, p. 150]). At the same time, the definition of parabolicity allows an equivalent characterization: a complete manifold $(M, g)$ is parabolic if any positive superharmonic function on $(M, g)$ is constant, otherwise, it is considered non-parabolic (see [24, p. 164]). Therefore, if a hypersurface $(M,g)$ is a complete parabolic Riemannian submanifold of $(\bar{M},\bar{g})$ and there exists a function $u(x) \in C^2(M)$ that has not zeros on $M$ and satisfying the inequality $\mathcal{L}\,u(x) \le 0$, then $u(x)$ is a constant and from (4.1) we deduce that $\varphi = 0$ and $\overline{Ric}(N_x, N_x) = 0$ at every point of $M$. In this case, $(M,g)$ must be a totally geodesic submanifold of $(\bar{M},\bar{g})$ and $\overline{Ric}(N_x, N_x)$ also vanishes at every point $x \in M$. This leads to the following straightforward theorem.

**Theorem 4.1.** *Let $(M,g)$ be a complete stable minimal hypersurface in $(n+1)$-dimensional Riemannian manifold $(\bar{M},\bar{g})$ with nonnegative Ricci curvature $\overline{Ric}$ in the direction of the unite normal vector $N_x$ to $M$ is non-negative at every point $x \in M$. If $(M,g)$ is parabolic (in particular, a complete manifold of finite volume) and there exists a function $u(x) \in C^2(M)$ that has no zeros on $M$ and satisfying the inequality $\mathcal{L}\,u(x) \le 0$ at each point of $M$, where $\mathcal{L}$ denotes the stability (or Jacobi) operator, then $(M,g)$ must be totally geodesic submanifold. Moreover, the Ricci curvature $\overline{Ric}(N_x, N_x)$ also vanishes at every point $x \in M$.*

There is a generalization of parabolic manifolds known as *stochastically complete manifolds*. Specifically, any parabolic manifold is stochastically complete, but the converse is not necessarily true. Recall that if we consider a minimal Wiener process on $(M,g)$, that is, a diffusion process generated by the Laplace–Beltrami operator $\Delta$ with absorption conditions at infinity, then if the probability of absorption at infinity in a finite amount of time is zero, the manifold $(M,g)$ is said to be stochastically complete (see, for instance, [33]). In conclusion, we can formulate the theorem.

**Theorem 4.2.** *Let $(M,g)$ be a complete, stable, and stochastically complete minimal hypersurface in a Riemannian manifold $(\bar{M},\bar{g})$, such that the Ricci curvature of $(\bar{M},\bar{g})$ in the direction of the unite normal vector $N_x$ to $M$ is non-negative at every point $x \in M$. If there exists a function $u(x) \in C^2(M) \cap L^1(M)$ that it has no zeros and satisfies the inequality $\mathcal{L}\,u(x) \le 0$ at each point of $M$, where $\mathcal{L}$ denotes the stability (or Jacobi) operator, then $(M,g)$ is totally geodesic submanifold.*

*Proof.* Let $(M, g)$ be a complete, stable, and stochastically complete minimal hypersurface in a Riemannian manifold $(\bar{M}, \bar{g})$, such that the Ricci curvature of $(\bar{M}, \bar{g})$ in the direction of the unite normal vector $N_x$ to $M$ is non-negative at every point $x \in M$.

Next, we recall that if $(M, g)$ a stochastically complete manifold, then every positive superharmonic function $u(x) \in L^1(M)$ is constant (see [33]). Therefore, if a hypersurface $(M, g)$ in $(\bar{M}, \bar{g})$ is a stochastically complete Riemannian manifold, then from (4.2) for a function $u(x) \in C^2(M) \cap L^1(M)$ that has no zeros on $M$, we obtain $\varphi = 0$. Consequently, $(M, g)$ must be a totally geodesic submanifold of $(\bar{M}, \bar{g})$.

**Remark 4.1.** In [10] and [34], with references to [35] and [36], respectively, it was claimed that on a complete, non-compact stable minimal hypersurface $(M, g)$ of $(\bar{M}, \bar{g})$, there exists a globally defined positive function $u(x)$ that belongs to the kernel of the Jacobi operator $\mathcal{L}$. In this case, if the Ricci curvature of $(\bar{M}, \bar{g})$ in the direction of the unite normal vector $N_x$ to $M$ is non-negative at every point $x \in M$, then $(\Delta u)(x) = -\left(\|\varphi\|_g^2(x) + \overline{Ric}(N_x, N_x)\right)u(x) \leq 0$. Therefore, $u(x)$ is a positive superharmonic function. This allows us to simplify the requirements of Theorems 4.1 and 4.2.